\documentstyle[12pt]{article}

\begin{document}
\newcommand{\ol }{\overline}
\newcommand{\ul }{\underline }
\newcommand{\ra }{\rightarrow }
\newcommand{\lra }{\longrightarrow }
\newcommand{\ga }{\gamma }
\newcommand{\st }{\stackrel }
\newcommand{\scr }{\scriptsize }
\title{\Large\textbf{Some Inequalities for Nilpotent Multipliers of Powerful $p$-Groups}}
\author{\textbf{Behrooz Mashayekhy\footnote{Corresponding author: mashaf@math.um.ac.ir} and Fahimeh Mohammadzadeh} \\
Department of Mathematics,\\ Center of Excellence in Analysis on Algebraic Stuctures,\\
Ferdowsi University of Mashhad,\\
P. O. Box 1159-91775, Mashhad, Iran.}
\date{ }
\maketitle
\begin{abstract}
In this paper we present some inequalities for the order, the
exponent, and the number of generators of the c-nilpotent multiplier
(the Baer invariant with respect to the variety of nilpotent groups
of class at most $c \geq 1$) of a powerful p-group. Our results
extend some of Lubotzky and Mann$^{,}$s (Journal of Algebra, 105
(1987), 484-505.) to nilpotent multipliers. Also, we give some
explicit examples showing the tightness of our results and
improvement some of the previous inequalities.
\end{abstract}
A.M.S. Classification 2000: 20C25, 20D15, 20E10, 20F12.\\
\textit{Keywords}: Nilpotent multiplier, Powerful $p$-group.
\newpage
\begin{center}
\hspace{-0.65cm}\textbf{1. Introduction and Motivation}\\
\end{center}

Let $G$ be a group with a free presentation $F/R$. The abelian
group
$$M^{(c)}(G)=\frac{R \cap \gamma_{c+1}(F)}{[R,\!\!\ _cF]}$$ is said to
be the $c$-nilpotent multiplier of $G$ (the Baer invariant of $G$,
after R. Baer [1], with respect to the variety of nilpotent groups
of class at most $c\geq 1$). The group $M(G)=M^{(1)}(G)$ is more
known as the Schur multiplier of $G$. When $G$ is finite, $M(G)$
 is isomorphic to the second cohomology group $H^{2}(G,\mathbf{C}^*)$
 [8].

It was conjectured for some time that the exponent of the Schur
multiplier of a finite $p$-group is a divisor of the exponent of the
group itself. I. D. Macdonald, J. W. Wamsley, and others [2] have
constructed an example of a group of exponent 4, whereas its Schur
multiplier has exponent 8, hence the conjecture is not true in
general. In 2007 P. Moravec [15] proved that if $G$ is a group of
exponent 4, then $\exp(M(G))$ divides 8. In 1973 Jones [7] proved
that the exponent of the Schur multiplier of a finite $p$-group of
class $c \geq 2$ and exponent $p^{e}$ is at most $p^{e(c-1)}$. A
result of G. Ellis [4] shows that if $G$ is a $p$-group of class $k
\geq 2$ and exponent $p^{e}$, then $\exp(M^{(c)}(G))\leq p^{e\lceil
k/2\rceil}$, where $\lceil k/2\rceil$ denotes the smallest integer
$n$ such that $n \geq k/2$. For $c=1$ P. Moravec [15] showed that
$\lceil k/2\rceil$ can be replaced by $2\lfloor \log_2{k}\rfloor$
which is an improvement if $k\geq 11$. Also he proved that if $G$ is
a metabelian group of exponent $p$, then $\exp(M(G))$ divides $p$.
S. Kayvanfar and M.A. Sanati [9] proved that $\exp(M(G))\leq
\exp(G)$ when $G$ is a finite $p$-group of class 3, 4 or 5 under
some arithmetical conditions on $p$ and the exponent of $G$. On the
other hand, the authors in a joint paper [13] proved that if $G$ is
a finite $p$-group of class $k$ with $p>k$, then
$\exp(M^{(c)}(G))|\exp(G)$. In 1972 Jones [6] showed that the order
of the Schur multiplier of a finite $p$-group of order $p^{n}$ with
center of exponent $p^{k}$ is bounded by $p^{(n-k)(n+k-1)/2}$. In
particular, $|G'||M(G)|\leq p^{\frac{n(n-1)}{2}}$. In 1973 Jones [7]
gave a bound for the number of generators of the Schur multiplier of
a finite $p$-group of class $c$ and special rank $r$. Recently the
authors in a joint paper [13] have extended this result to the
$c$-nilpotent multipliers. In 1987 Lubotzky and Mann [10] presented
some inequalities for the Schur multiplier of a powerful $p$-group.
They gave a bound for the order, the exponent and the number of
generators of the Schur multiplier of a powerful $p$-group. Their
results improve the previous inequalities for powerful $p$-groups.
In this paper we will extend some results of Lubotzky and Mann [10]
to the nilpotent multipliers and give some upper bounds for the
order, the exponent and the number of generators of the
$c$-nilpotent multiplier of a $d$-generator powerful $p$-group $G$
as follows:
 $$d(M^{(c)}(G))\leq \chi_{c+1}(d),\ \exp(M^{(c)}(G))|
\exp(G),$$ $$ {\rm and}\  |M^{(c)}(G)|\leq
p^{\chi_{c+1}(d)}\exp(G),$$
 where $\chi _{c+1}
(d)$ is the number of basic commutators of weight $c+1$ on $d$
letters [5]. Our method is similar to that of [10]. Finally, by
giving some examples of groups and computing the number of
generators, the order and the exponent of their $c$-nilpotent
multipliers explicitly, we compare these numbers with the bounds
obtained and show that our results improve some of the previously mentioned inequalities.\\

\begin{center}
\textbf{2. Notation and Preliminaries}
\end{center}

Here we will give some definitions and theorems that will be used in
our work. Throughout this paper $\mho_i(G)$ denotes the subgroup of
$G$ generated by all $p^{i}$th powers, $P_i(G)$ is defined by:
$P_1(G)=G$, and $P_{i+1}(G)=[P_i(G),G]\mho_1(P_i(G))$. Finally
$d(G)$, $cl(G)$, $l(G)$, $sr(G)$ denote respectively, the minimal
number of generators, the nilpotency class, the derived length and
the special rank of $G$, while $e(G)$ is
defined by $\exp(G)=p^{e(G)}$.\\
\textbf{\bf Theorem 2.1} (M. Hall [5]). \textit{Let $F$ be a free
group on $\{x_{1}, x_{2}, ...,x_{d}\}$. Then for all $1\leq i \leq
n$, $$ \frac{\gamma_{n}(F)}{\gamma_{n+i}(F)} $$ is a free abelian
group freely generated by the basic commutators of weights $ n,
n+1,..., n+i-1 $  on the letters $\{x_{1}, x_{2}, ...,x_{d}\}$
(for a definition of basic commutators see [5])}.\\
\textbf{\bf Lemma 2.2} (R. R. Struik [16]).\textit{ Let $\alpha$ be
a fixed integer and $G$ be a nilpotent group of class at most $n$.
If $b_j \in G$ and $r<n$, then
$$[b_1,..,b_{i-1},b_i^{\alpha},b_{i+1},...,b_r]=[b_1,...,b_r]^{\alpha}c_1^{f_1(\alpha)}c_2^{f_2(\alpha)}...,$$
where the $c_k$ are commutators in $b_1,...,b_r$ of weight strictly
greater than $r$, and every $b_j$, $1\leq j \leq r$, appears in each
commutator $c_k$, the $c_k$ listed in ascending
order. The $f_i$ are of the following form: \\
$$f_i(n)=a_1{n \choose 1}+ a_2{n \choose
2}+...+a_{w_i}{n \choose w_i},$$ with $a_j \in \bf{Z}$, and $w_i$ is
the weight of $c_i$ ( in the $b_i$ ) minus $(r-1)$.}

Powerful $p$-groups were formally introduced in [10]. They have
played a role in the proofs of many important results in $p$-groups.
We will discuss some of them in this section. A $p$-group $G$ is
called \emph{powerful} if $p$ is odd and $G' \leq \mho_1(G)$ or
$p=2$, and $G' \leq \mho_2(G)$. There is a related notion that is
often used to find properties of powerful $p$-groups. If $G$ is a
$p$-group and $H \leq G$, then $H$ is said to be \emph{powerfully
embedded} in $G$ if $[G,H]\leq \mho_1(H)$ ($[G,H]\leq \mho_2(H) $
for $p=2$). Any powerfully embedded subgroup is itself a powerful
$p$-group and must be normal in the whole group. Also a $p$-group is
powerful exactly when it is powerfully embedded in itself. While it
is obvious that factor groups and direct products of powerful
$p$-groups are powerful, this property is not subgroup-inherited
[10].

We will require some standard properties of powerful $p$-groups.
For the sake of convenience we collect them here.\\
\textbf{\bf Theorem 2.3} ([10]).\textit{ The following statements
hold for a powerful $p$-group
$G$.\\
i) $\gamma_i(G), G^{i}, \mho_i(G), \Phi(G)$ are powerfully
embedded in $G$.\\
ii) $P_{i+1}(G)=\mho_i(G)$ and
$\mho_i(\mho_j(G))=\mho_{i+j}(G)$.\\
iii) Each element of $\mho_i(G)$ can be written as $a^{p^{i}},$
for
some $a\in G$ and hence $\mho_i(G)=\{g^{p^{i}}: g\in G\}$.\\
iv) If $G=\langle a_1,a_2,...,a_d\rangle$, then $\mho_i(G)=\langle
a_1^{p^{i}},a_2^{p^{i}},...,a_d^{p^{i}}\rangle$.\\
v) If $H \subseteq G$, then $d(H) \leq d(G)$}.\\
\textbf{\bf Proposition 2.4} ([10]). \textit{Let $N$ be a powerfully
embedded subgroup of $G$. If $N$ is the normal closure of
some subset of $G$, then $N$ is actually generated by this subset.}\\
\textbf{\bf Lemma 2.5.} \textit{Let $H, K$ be normal subgroups of
$G$ and $H \leq K[H,G]$. Then $H \leq K[H,\  _l \ G]$ for any $l\geq
1$. In particular, if $G$ is nilpotent, then $H \leq K$ .}\\
\textit{Proof.} An easy exercise.\ $\Box$\\
\textbf{\bf Lemma 2.6.} \textit{Let $G$ be a finite $p$-group and $N
\unlhd G$. Then $N$ is powerfully
 embedded in $G$ if and only if $N/[N,G,G]$ is powerfully embedded
 in $G/[N,G,G]$.}\\
 \textit{Proof.} See a remark in the proof of Theorem 1.1 in
 [10].\ $\Box$\\
\textbf{\bf Remark 2.7.} To prove that a normal
subgroup $N$ is powerfully embedded in $G$ we can assume that\\
i) $[N,G,G]=1$ by the above lemma.\\
ii) $\mho_1(N)=1$ ( $\mho_2(N)=1$ for $p=2$ ) and try to show that
$[N,G]=1$.\\
iii) $[N,G]^{2}=1$ whenever $p=2$, since if we assume that $
N/[N,G]^{2}$ is powerfully embedded in  $ G/[N,G]^{2}$, then $N$ is
powerfully embedded in $G$. This follows from the proof of Theorem 4.1.1 in [10].\\
\newpage
\begin{center}
\textbf{3. Main Results}
\end{center}

In order to prove the main results we need the following theorem.\\
\textbf{\bf Theorem 3.1}. \textit{Let $F/R$ be a free presentation
of a powerful $d$-generator $p$-group $G$. Let $Z=R/[R,\!\!\ _cF]$
and $H=F/[R,\!\!\ _cF]$, so that $G\cong H/Z$. Then $\gamma
_{c+1}(H)$ is powerfully embedded in $H$ and $d(\gamma
_{c+1}(H))\leq \chi_{c+1}(d)$.}\\
\textit{Proof.} First let $p$ an odd prime. We may assume that
$\mho_1 (\gamma _{c+1}(H))=1$ and try to show that $[(\gamma
_{c+1}(H)),H]=1$ by Remark 2.7(ii). Also we may assume that $\gamma
_{c+3}(H)=1$ by Remark 2.7(i).  Let $a,b_1,b_2,...,b_c\in H$. Then
by Lemma 2.2,
$$[a^{p},b_1,...,b_c]=[a,b_1,...,b_c]^{p}c_1^{f_1(p)}c_2^{f_2(p)}...\ \ \ .$$
Since $\gamma_{c+3}(H)=1$ and $\mho_1 (\gamma _{c+1}(H))=1$ we have
$[a,b_1,...,b_c] ^{p}=1, c_i^{f_i(p)}=1 $, for all $i \geq 2$. Also
$p>2$ implies that $p| f_1(p)$, and hence $c_1^{f_1(p)}=1$ so
$a^{p}\in Z_c(H)$ and $\mho_1(H)\subseteq Z_c(H)$. The powerfulness
of $G$ yields $H' \leq \mho_1(H)Z \leq Z_c(H)$. Therefore $[H',\!\!\
_cH]=1$, as desired. Since $H/Z$ is generated by $d$ elements and
$Z\leq Z_c(H)$, $\gamma _{c+1}(H)$ is the normal closure of the
commutators  of weight $c+1$ on $d$ elements. Hence Proposition 2.4
completes the proof, for $p>2$.

If $p=2$, then the proof is similar, so we leave out the details,
but note that in this case
$$[a^{4},b_1,...,b_c]=[a,b_1,...,b_c]
^{4}c_1^{f_1(4)}c_2^{f_2(4)}...\ \ \ .$$ By Remark 2.7 we can assume
$\gamma_{c+3}(H)=\mho_2(\gamma_{c+1}(H)) =([\gamma_{c+1}(H),H])
^{2}=1$. Hence we have  $[a^{4},b_1,...,b_c]=1 \ (c_1^{f_1(4)}=1$,
since $2|f_1(4)$) so $\mho_2(H)\subseteq Z_c(H). \ \Box$

An interesting corollary of this theorem is as follows.\\
\textbf{\bf Corollary 3.2.} \textit{Let $G$ be powerful $p$-group
with $d(G)=d$. Then
$d(M^{(c)}(G))\leq \chi_{c+1}(d)$.} \\
\textit{Proof.} Let $F/R$ be a free presentation of $G$ with
$Z=R/[R,\!\!\!\ _cF]$, so that $G\cong H/Z$, where $H=F/[R,\!\!\
_cF]$. Then the above result and Theorem 2.3(v) implies that
$$d(\frac{R \cap \gamma_{c+1}(F)}{[R,\!\!\ _cF]})\leq d(\frac{
\gamma_{c+1}(F)}{[R,\!\!\ _cF]})\leq \chi_{c+1}(d).$$ Hence the
result follows. $\Box$

Note that by a similar method we can prove Corollary 2.2 of [10]
without using the concept of covering group for $G$.

The authors in a joint paper [12] have proved that if G is a finite
d-generator p-group of special rank r and nilpotency class $t$, then
$ d(M^{(c)} (G)) \leq \chi _{c+1} (d) + r ^{c+1} (t-1) $. Clearly
Corollary 3.2
improves this bound for nonabelian powerful $p$-groups.\\
\textbf{\bf Theorem 3.3.} \textit{ Let $G$ be powerful $p$-group.
Then $e(M^{(c)}(G))\leq e(G)$.} \\
\textit{Proof.} Let $p>2$ and $F/R$ be a free presentation of $G$
with $Z=R/[R,\!\!\ _cF]$ and $H=F/[R,\!\!\ _cF]$, so that $G\cong
H/Z$. Since $e(R \cap \gamma_{c+1}(F)/[R,\!\!\ _cF])\leq e(
\gamma_{c+1}(H))$ and $e(H/Z_c(H))\leq e(G)$ it is enough to show
that $e( \gamma_{c+1}(H))=e(H/Z_c(H))$. We will establish by
induction on $k$ the equality
\begin{equation}
\mho_k(\gamma _{c+1}(H))=[\mho_k(H),\!\!\ _cH],
\end{equation}
which implies the above claim.

If $k=0$, then (1) holds. Now Assume that (1) holds for some $k$.
Since $\gamma_{c+1}(H)$ is powerfully embedded in $H$ by Theorem
3.1, we have $\mho_{k+1}(\gamma _{c+1}(H))=\mho_1 (\mho_k(\gamma
_{c+1}(H)))$, by Theorem 2.3(ii). Similarly $\mho_{k+1}(G))=\mho_1
(\mho_k(G)))$. Since $G\cong H/Z $ we have $\mho_{k+1}(H)Z/Z=\mho_1
(\mho_k(H)Z)Z/Z$. Therefore $$[\mho_{k+1}(H),\!\!\ _cH]=
[\mho_{k+1}(H)Z,\!\!\ _cH]=[\mho_1 (\mho_k(H)Z)Z,\!\!\  _cH]$$ $$=
[\mho_1 (\mho_k(H)Z),\!\!\ _cH].$$ This implies that
\begin{equation}
[\mho_{k+1}(H),\!\!\ _cH]=[\mho_1 (\mho_k(H)Z),\!\!\ _cH].
\end{equation}
Thus (1) for $k+1$ is equivalent to $\mho_1 (\mho_k(\gamma
_{c+1}(H)))=[\mho_1 (\mho_k(H)Z),\!\!\ _cH]$. Since $\mho_k(\gamma
_{c+1}(H))$ is powerfully embedded in $H$ by Theorem 2.3(i), this
implies, by (1) and Lemma 2.2,
\begin{eqnarray*}
[\mho_1 (\mho_k(H)Z),\!\!\ _cH]&\leq& \mho_1 ([\mho_k(H)Z,\!\!\
_cH])[\mho_k(H)Z,\!\!\ _cH,H]\\&\leq& \mho_1 ([\mho_k(H),\!\!\
_cH])[\mho_k(H),\!\!\ _cH,H]\\&\leq& \mho_1 (\mho_k(\gamma
_{c+1}(H)))[\mho_k(\gamma_{c+1}(H)),H]\\&\leq& \mho_1 (\mho_k(\gamma
_{c+1}(H))).
\end{eqnarray*}

For the reverse inclusion note that since $\mho_1 (\mho_k(\gamma
_{c+1}(H)))=\mho_1 ([\mho_k(H), \ _cH])$ it is enough to show that
$$\mho_1 ([\mho_k(H),\!\!\ _cH]) \equiv 1 \pmod{[\mho_1 ( \mho_k(H)Z),\!\!\ _cH]}.$$
By Theorem 2.3(i), $\mho_k(H/Z)$ is powerfully embedded in  $H/Z$ so
that
\begin{equation}
[\frac{\mho_{k}(H)Z}{ Z}, \frac{H}{Z}]\leq  \frac{\mho_1(
\mho_{k}(H)Z)Z}{Z}.
\end{equation}
Also (2) implies that $\mho_1(\mho_{k}(H)Z)\leq Z_c(H)
\pmod{[\mho_{k+1}(H),\ _cH]}.$ Now (2), (3) and the last inequality
imply that
$$[\mho_{k}(H)Z,H]\leq \mho_1(\mho_{k}(H)Z)Z \leq Z_c(H)
\pmod{[\mho_{k+1}(H), \ _cH]}.$$ Hence by Lemma 2.2
\begin{eqnarray*}
\mho_1([\mho_{k}(H),\ _cH])&\equiv& \mho_1([\mho_{k}(H)Z,\!\!\
_cH])\\&\equiv& [\mho_1 (\mho_{k}(H)Z),\!\!\ _cH]
\\&\equiv& 1 \pmod{[\mho_1 (\mho_{k}(H)Z),\!\!\ _cH]},
\end{eqnarray*}
as desired.

If $p=2$, then the proof is similar to the previous case. This
completes the proof.\ $\Box$

Note that G. Ellis [4], using the nonabelian tensor products of
groups, showed that $\exp(M^{(c)}(G))$ divides $\exp(G)$ for all
$c\geq 1$ and all $p$-groups satisfying $[[G^{p{i-1}},G],G]\subseteq
G^{p^i}$ for $1\leq i\leq e$, where $\exp(G)=p^e$. Note that the
results of [10] imply that every powerful $p$-group $G$ satisfies
the latter commutator condition.

Lubotzky and Mann [10] found  bounds for $cl(G), l(G), |G|$ and
$|M(G)|$ of a powerful $d$-generator $p$-group $G$  of  exponent $p
^{e}$ as follows:
$$cl(G)\leq e,\ l(G)\leq \log_2{e}+1,\ |G|\leq p^{de}\ {\rm and} \  |M (G)|\leq p^{
(d(d-1)/2)e}.$$ In the following proposition we find an upper bound
for the
 order of $c$-nilpotent multiplier of $G$.\\
\textbf{\bf Proposition 3.4.} \textit{Let $G$ be a powerful
$p$-group, with $d(G)=d $ and $e(G)=e$. Then $|M^{(c)}(G)|\leq
p^{\chi_{c+1}(d)e}$}.\\
\textit{Proof.} It is obtained by combining Corollary 3.2 and
Theorem 3.3.\ $\Box$

\begin{center}
\textbf{4. Some Examples}
\end{center}

 In this final section we are going to give some explicit examples
 of $p$-groups and calculate their $c$-nilpotent multipliers in
 order to compare our new bounds with the exact values. This will show
 tightness of our results and improvement some of the previously
 mentioned inequalities.\\
\textbf{Example 4.1.} Let $G$ be a finite abelian $p$-group.
 Clearly $G$ is a powerful $p$-group and by the fundamental theorem
 of finitely generated abelian groups $G$ has the following
 structure
 $$ G\cong {\mathbf Z}_{p^{\alpha_1}}\oplus {\bf Z}_{p^{\alpha_2}}\oplus \ldots \oplus {\bf Z}_{p^{\alpha_d}} $$
for some positive integers $\alpha_1, \alpha_2, \ldots , \alpha_d$,
where $\alpha_1\geq \alpha_2\geq \ldots \geq \alpha_d$. By [11] the
$c$-nilpotent multiplier of $G$ can be calculated explicitly as
follows:
$$ M^{(c)}(G)\cong {\bf Z}_{p^{\alpha_2}}^{(b_2)}\oplus {\bf
Z}_{p^{{\alpha_3}}}^{(b_3-b_2)}\oplus \ldots \oplus {\bf
Z}_{p^{\alpha_d}}^{(b_d-b_{d-1})},$$ where $b_i=\chi_{c+1}(i)$ and
${\bf Z}_n^{(m)}$ denotes the direct sum of $m$ copies of the cyclic
group ${\bf Z}_n$. Now it is easy to see that\\
$(i)$ $d(M^{(c)}(G))=\chi_{c+1}(d)$, where $d=d(G)$. Hence the bound
of Corollary 3.2 is attained and the best one in the abelian case.\\
$(ii)$ $e(M^{(c)}(G))=\alpha_2$, whereas $e(G)=\alpha_1$. Hence the
bound of Theorem 3.3 is attained when $\alpha_1=\alpha_2$ and it is
the best one in the abelian case.\\
$(iii)$
$|M^{(c)}(G)|=p^{\alpha_2b_2+\sum_{i=3}^{d}\alpha_i(b_i-b_{i-1})}\leq
p^{\alpha_1}\chi_{c+1}(d)$. Hence the bound of Proposition 3.4 is
attained if and only if $\alpha_1=\alpha_2=\ldots =\alpha_d$.\\
\textbf{Example 4.2.} Let $p$ be any odd prime number and $s,\ t$ be
positive integers with $s\geq t$. Consider the following finite
$d$-generator $p$-group with nilpotency class $2$:
$$P_{s,t}=\langle y_1,\ldots ,
y_d:y_i^{p^s}=[y_j,y_k]^{p^t}=[[y_j,y_k],y_i]=1,\ 1\leq i,j,k\leq
d,\ j\neq k\rangle.$$ One can see that $P_{s,t}$ is not a powerful
$p$-group (clearly $\mho_1(P_{1,1})=1$). By [14] the $c$-nilpotent
multiplier of $P_{s,t}$ is as follows:
$$ M^{(c)}(P_{s,t})\cong {\bf Z}_{p^s}^{(\chi_{c+1}(d))}\oplus {\bf
Z}_{p^t}^{(\chi_{c+2}(d))}.$$ Therefore we have\\
$(i)$
$d(M^{(c)}(P_{s,t}))=\chi_{c+1}(d)+\chi_{c+2}(d)>\chi_{c+1}(d)$.
Hence the condition of being powerful cannot be omitted from
Corollary 3.2.\\
$(ii)$
$|M^{(c)}(P_{s,t}))|=p^{s\chi_{c+1}(d)+t\chi_{c+2}(d)}>p^{s\chi_{c+1}(d)}$.
Hence powerfulness is also a necessary condition for the bound of
Proposition 3.4. Note that here we have
$e(M^{(c)}(P_{s,t}))=se(P_{s,t})$.

The authors in a joint paper [13] have proved that
$\exp(M^{(c)}(G))|\exp(G)$, when $G$ is a nilpotent $p$-group of
class $k$, and $k<p$. In the following example we find a powerful
$p$-group of class $k \geq p$ such that $\exp(M^{(c)}(G))$ divides
$\exp(G)$.\\
\textbf{\bf Example 4.3} ([17]). We work in ${\rm GL}({\bf Z}
_{p^{l+2}})$, the $2\times2$ invertible matrices over the ring of
integers modulo $p^{l+2}$. In this ring any integer not divisible
by $p$ is invertible. Consider the matrices \\
$$X=\left[\begin{array}{ll}1&0\\
0&1-p
\end{array}\right],Y=\left[\begin{array}{ll}1/(1-p)&p/(1-p)\\
0&1
\end{array}\right],Z=\left[\begin{array}{ll}1&p\\
0&1
\end{array}\right].$$
One quickly calculates that $[X,Y]=Z^{p}$, $[X,Z]=Z^{p}$,
$[Y,Z]=Z^{p}$ and
\begin{equation}
[Z^{p},\!\!\ _kX]=\left[\begin{array}{ll}1&(-1)^{k+2}p^{k+2}\\
0&1
\end{array}
\right].\end{equation} Notice also that
$X^{p^{l+1}}=Y^{p^{l+1}}=Z^{p^{l+1}}=1$. We claim that $P=\langle
X,Y,Z\rangle$ is a powerful $p$-group. By the above relations we can
express every word in $P$ as a product $X^{a}Y^{b}Z^{c}$ for some $
0\leq a,b,c < p^{l+1}$. Also
$$X^{a}Y^{b}Z^{c}=\left[\begin{array}{ll}\frac{1}{(1-p)^{b}}&\frac{1+pc-(1-p)^{b}}{(1-p)^{b}}\\
0&(1-p)^{a}\end{array}\right]$$ and hence all of these elements are
distinct. Therefore the order of $P$ is $ p^{3(l+1)}$ and hence $P$
is a $p$-group and the relations imply that $P'\leq \mho_1(P)$.
Therefore $P$ is a powerful $p$-group. The exponent of $P$ is
$p^{l+1}$, and (4) implies that $P$ has nilpotency class $l+1$. By
Theorem 3.3 $\exp(M^{(c)}(P))$ divides $\exp(P)$. Note that the
nilpotency class of $P$ is $l+1$ which is greater than or equal to
$p$.

Let $G$ be a finite d-generator p-group of order $p^n$ where $p$ is
any prime. By [11] we have
 $$p^{\chi_{c+1}(d)}\leq |M^{(c)}(G)||\gamma_{c+1}(G)|\leq p^{\chi_{c+1}(n)}.$$
Now if we put $l=2$ in the above example, then $P$ is $3$-generator
powerful $p$-group of order $p^9$ with nilpotency class $3$. Thus by
the above bounds we have
$$ p^{18}=p^{\chi_4(3)}\leq |M^{(3)}(P)||\gamma_4(P)|=|M^{(3)}(P)|\leq
p^{\chi_4(9)}=p^{1620}.$$ But by Proposition 3.4 $|M^{(3)}(P)|\leq
p^{3\chi_4(3)}=p^{54}$. Hence this example and also Example 4.1 show
that Proposition 3.4 improves the above bound for powerful
$p$-groups.\\
\textbf{Example 4.4.} Using the list of nonabelian groups of order
at most $30$ with their $c$-nilpotent multipliers for $c=1,2$ in the
table of Fig.2 in [3], we are going to give two nonabelian powerful
$p$-groups in order to compute explicitly the number of generators,
the order and the exponent of their $2$-nilpotent multipliers and
then compare these numbers with bounds obtained.\\
$(i)$ Consider the finite $2$-group $G=\langle
a,b:a^2=1,aba=b^{-3}\rangle$. It is easy to see that $G$ is a
powerful $2$-group and $|G|=16$, $d(G)=2$, $\exp(G)=8$. By [3,
Fig.2, $\sharp$ 13] $M^{(2)}(G)\cong {\bf Z}_2^{(2)}$ and hence
$|M^{(2)}(G)|=4$, $d(M^{(2)}(G))=2$, $\exp(M^{(2)}(G))=2$. It is
seen that the bound of Corollary 3.2 is attained.\\
$(ii)$ Consider the finite $3$-group $G=\langle
a,b:a^3=1,a^{-1}ba=b^{-2}\rangle$. It is easy to see that $G$ is a
powerful $3$-group and $|G|=27$, $d(G)=2$, $\exp(G)=9$. By [3,
Fig.2, $\sharp$ 40] $M^{(2)}(G)\cong {\bf Z}_3^{(2)}$ and hence
$|M^{(2)}(G)|=9$, $d(M^{(2)}(G))=2$, $\exp(M^{(2)}(G))=3$. It is
also seen that the bound of Corollary 3.2 is attained.\\
\begin{center}
\textbf{Acknowledgements}
\end{center}
This research was in part supported by a grant from Center of
Excellence in Analysis on Algebraic Structures, Ferdowsi University
of Mashhad.\\
The authors are grateful to the referee for useful comments and
careful corrections.

\end{document}